\documentclass[12pt,twoside,a4paper]{amsart}
\usepackage{amssymb}
\usepackage{graphicx}

\date{\today}

\usepackage[latin1]{inputenc}
\usepackage[T1]{fontenc}

\def\dbar{\bar\partial}

\def\R{{\mathbb R}}
\def\C{{\mathbb C}}
\def\P{{\mathbb P}}

\def\O{{\mathcal O}}

\def\Re{{\rm Re\,  }}

\def\sgn{\text{sgn}\;}
\def\codim{\text{codim}\,}
\def\ann{\text{ann}\,}

\def\1{\mathbf 1}

\def\Z{{\mathbb Z}}
\def\m{{\mathfrak m}}
\def\J{{\mathcal J}}

\def\a{{\mathfrak a}}
\def\ord{\text {ord}}
\def\vE{\text {ord}_E}
\def\jac{\text {Jac}}
\def\I{{\mathcal I}}
\def\np{\text{NP}}
\def\1{\mathbf 1}

\def\res{\text{Res}}
\def\vol{\text{Vol}}

\def\be{\begin{equation}}
\def\ee{\end{equation}}

\newtheorem{thm}{Theorem}[section]
\newtheorem{lma}[thm]{Lemma}

\newtheorem{prop}[thm]{Proposition}

\newtheorem{question}[thm]{Question}

\newtheorem*{thmA}{Theorem A}
\newtheorem*{questionB}{Question B}
\newtheorem*{questionB'}{Question B'}
\newtheorem*{thmC}{Theorem C}

\theoremstyle{definition}

\theoremstyle{remark}
\newtheorem{preremark}[thm]{Remark}
\newtheorem{preex}[thm]{Example}

\newenvironment{ex}{\begin{preex}}{\qed\end{preex}}

\numberwithin{equation}{section}

\numberwithin{equation}{section}

\begin{document}

\title[On weighted Bochner-Martinelli residue currents]{On weighted Bochner-Martinelli \\ residue currents}

\date{\today}
\thanks{The author was partially supported by the Swedish Research Council and NSF grant DMS-0901073} 

\author{Elizabeth Wulcan}

\address{Dept of Mathematics, University of Michigan, Ann Arbor \\ MI 48109-1043\\ USA}

\email{wulcan@umich.edu}



\begin{abstract}
We study the weighted Bochner-Martinelli residue current $R^p(f)$ associated with a sequence $f=(f_1,\dots,f_m)$ of holomorphic germs at $0\in\mathbf{C}^n$, whose common zero set equals the origin, and $p=(p_1,\ldots, p_m)\in\mathbb N^n$. Our main results are a description of $R^p(f)$ and its annihilator ideal in terms of the Rees valuations of the ideal generated by $(f_1^{p_1},\ldots, f_m^{p_m})$ and an explicit description of $R^p(f)$ when $f$ is monomial. For a monomial sequence $f$ we show that $R^p(f)$ is independent of $p$ if and only if $f$ is a regular sequence. 
\end{abstract}

\maketitle
\section{Introduction}
Let $f=(f_1,\ldots, f_m)$ be sequence of germs of holomorphic functions at $0\in\C^n$, such that $V(f)=\{f_1=\ldots =f_m=0\}=\{0\}$. If $f$ is a \emph{regular sequence}, that is, $m=n$, then there is a canonical residue (current) associated with $f$ - the Grothendieck residue $\text{Res}(\frac{\bullet}{f_1\cdots f_m})$, see ~\cite{GH}, and its current avatar the \emph{Coleff-Herrera product}
$R_{CH}(f)=\dbar[1/f_1]\wedge\cdots\wedge\dbar[1/f_m]$ 
introduced in ~\cite{CH}. 
In ~\cite{PTY} Passare-Tsikh-Yger constructed residue currents based on the Bochner-Martinelli kernel as a natural generalization of the Coleff-Herrera product. This idea is further developed in ~\cite{BY}, where Berenstein-Yger introduced \emph{weighted Bochner-Martinelli residue currents}. 

Let $p=(p_1,\ldots, p_m)\in\mathbb N^m$ and let $f^p$ denote the sequence $(f_1^{p_1},\ldots, f_m^{p_m})$; here $\mathbb N$ denotes the natural numbers $1,2,\ldots$. For each ordered multi-index $\I=\{i_1,\ldots,i_n\}\subseteq \{1,\ldots, m\}$ let 
\begin{equation}\label{blanc}
R^p_\I(f)=\dbar |f^p|^{2\lambda}\wedge
c_n \sum_{\ell=1}^n(-1)^{\ell-1}
\frac
{\bar{f_{i_\ell}}|f_i|^{2(p_i-1)}\bigwedge_{q\neq \ell}'
\dbar(\bar{f_{i_q}} |f_{i_q}|^{2(p_{i_q}-1)})}
{|f^p|^{2n}}\bigg |_{\lambda=0},
\end{equation}
where $c_n=(-1)^{n(n-1)/2}(n-1)!$, 
$|f^p|^2=|f_1^{p_1}|^2+\ldots+|f_m^{p_m}|^2$, $\bigwedge '$ denotes increasing order in $q$ in the wedge product, and $\alpha |_{\lambda =0}$ denotes the analytic continuation of the form $\alpha$ to $\lambda=0$. 
Moreover, let $R^p(f)$ denote the vector-valued current with entries $R^p_\I(f)$; we will refer to this as the \emph{Bochner-Martinelli residue current of weight $p$} associated with $f$. Then $R^p(f)$ is a well-defined $(0,n)$-current with support at the origin and $\overline g R^p_\I(f)=0$ if $g$ is a holomorphic function that vanishes at the origin. It follows that the coefficients of the $R^p_\I(f)$ are just finite sums of holomorphic derivatives at the origin. If $p=(1,\ldots, 1)$, then $R^p(f)$ is the \emph{Bochner-Martinelli residue current} associated with $f$, ~\cite{PTY}; we denote it by $R(f)$ and its entries by $R_\I(f)$. Note that, in fact, 
\begin{equation}\label{enkel}
R^p_\I(f)=f_{i_1}^{p_{i_1}-1}\cdots f_{i_n}^{p_{i_n}-1} R_\I(f^p).
\end{equation}
Indeed, the sequence $f^p$ in the factor $\dbar |f^p|^{2\lambda}$ in \eqref{blanc} can be replaced by any sequence of functions that vanish at the origin.

Let $\O^n_0$ be the local ring of germs of holomorphic functions at $0\in\C^n$. Given a germ of a current $\mu$ at $0\in\C^n$, let $\ann \mu$ denote the (holomorphic) \emph{annihilator ideal} of $\mu$, that is $\ann \mu=\{h\in\O^n_0, h\mu=0\}$. 
Our first result concerns $\ann R^p(f)$. 
 Let $\a (f)$ denote the ideal generated by the $f_i$ in $\O^n_0$. Recall that $h\in\O^n_0$ is in the \emph{integral closure} of $\a (f)$, denoted by $\overline{\a (f)}$, if $|h|\leq C|f|$, for some constant $C$. Moreover, recall that $\a (f)$ is a \emph{complete intersection ideal} if it can be generated by $n=\codim V(f)$ functions. Note that this condition is slightly weaker than that $f$ is a regular sequence. 
Also, recall that, given ideals $\a, \mathfrak b\subseteq \O_0^n$, the colon ideal $\a\colon \mathfrak b$ is the ideal $\a\colon \mathfrak b=\{h\in \O_0^n: h\mathfrak b\subseteq \a\}$.

We also provide a characterization of the non-vanishing entries of $R^p(f)$. 
Let $\pi:X\to (\C^n,0)$ be a log-resolution of $\a (f)$, see ~\cite[Def.~ 9.1.12]{Laz}. Following ~\cite{JW} say that a multi-index $\I=\{i_1,\ldots,i_n\}$ is \emph{essential} with respect to $f$ if there is an exceptional prime $E\subseteq\pi^{-1}(0)$ of $X$ such that the mapping $[f_{i_1}\circ \pi:\ldots :f_{i_n}\circ\pi]: E\to \C\P^{n-1}$ is surjective and moreover $\vE(f_{i_k})\leq \vE(f_{\ell})$ for $1\leq k\leq n, 1\leq\ell\leq m$, see Section ~\ref{background} and also ~\cite[Section~3]{JW} for details. The valuations $\vE$ are precisely the \emph{Rees valuations} of $\a (f)$. We say that $\I$ is $p$-\emph{essential} if it is essential with respect to $f^p$. For $h\in\O_0^n$, let $(h)$ denote the ideal generated by $h$.

\begin{thmA}
Suppose that $f$ is a sequence of germs of holomorphic functions at $0\in\C^n$, such that $V(f)=\{0\}$. Let $R^p(f)$ be the corresponding Bochner-Martinelli residue current of weight $p$. Then the entry $R^p_\I(f)\not\equiv 0$ if and only if $\I$ is $p$-essential. Moreover
\begin{equation}\label{eqa}
\bigcap_{\I ~~~~~~ p-\text{essential}}
\overline{\a (f^{p})^n}\colon (f_{i_1}^{p_{i_1}-1}\cdots f_{i_n}^{p_{i_n}-1})
\subseteq\ann R^p(f)\subseteq \a (f).
\end{equation}
The left inclusion in \eqref{eqa} is strict whenever $n\geq 2$. If the right inclusion is an equality, then $\a (f)$ is a complete intersection ideal. 
\end{thmA}
The new results in Theorem ~A are the characterization of the non-vanishing entries and the last two statements. Berenstein-Yger ~\cite{BY} showed that $\overline{\a (f_j^{p_j})^n}\colon (f_{i_1}^{p_{i_1}}\cdots f_{i_n}^{p_{i_n}})
\subseteq\ann R^p_\I(f)$, and it is easy to see from Andersson's construction of residue currents in ~\cite{A} that the right inclusion in \eqref{eqa} holds.
In fact, Berenstein-Yger defined currents $R^p_\I(f)$ also when $\dim V(f)>0$. The inclusions \eqref{eqa} hold true also in this case, and one can even replace the leftmost ideal by $\bigcap_{\I ~~~~~~ p-\text{essential}}
\overline{\a (f_j^{p_j})^\mu}\colon (f_{i_1}^{p_{i_1}}\cdots f_{i_n}^{p_{i_\mu}})$, where $\mu=\min(m,n)$.

 Also, for $R(f)=R^{(1,\ldots, 1)}(f)$ Theorem ~A was proved in parts in ~\cite{PTY}, ~\cite{A}, and ~\cite{JW}. If $f$ is a regular sequence, then the only entry $R_{\{1,\ldots, m\}}(f)$ of $R(f)$ coincides with the Coleff-Herrera product $R_{CH}(f)$, whose annihilator ideal is precisely $\a (f)$, see ~\cite{DS, P}. This should be compared to ~\cite[Chapter~5.1]{GH} where $\res(\frac{\bullet}{f_1\cdots f_m})$ is defined using the Bochner-Martinelli kernel. 
The idea of regarding (complete intersection) ideals of holomorphic functions as the annihilator ideals of certain residue currents is central for many applications, see ~\cite{BGVY}. For $p=(1,\ldots, 1)$,  the inclusions \eqref{eqa} read $\overline{\a (f)^n}\subseteq\ann R(f)\subseteq \a (f)$, which gives a direct proof of the Brian{\c c}on-Skoda Theorem ~\cite{BS}: $\overline{\a (f)^n}\subseteq \a (f)$. For other applications of Bochner-Martinelli residue currents, see for example ~\cite{AG}, ~\cite{ASS}, and ~\cite{VY}.

Weighted Bochner-Martinelli residue currents were introduced in ~\cite{BY} as a tool to construct Green currents but also as a natural extension of Bochner-Martinelli residue currents in the spirit of Lipman ~\cite{Lipman}; the currents have been further studied in ~\cite{BVY} and ~\cite{Y}. In the monograph ~\cite{Lipman} 
not only the residue $\res(\frac{\bullet}{f_1\cdots f_m})$ associated with a sequence $f$ plays a role but also residues of the form $\res(\frac{f_1^{p_1-1}\cdots f_m^{p_m-1}\bullet}{f_1^{p_1}\cdots f_m^{p_m}})$. 
The currents $R^p(f)$ can thus be seen as analogues of this list of residues. If $f$ is a regular sequence, then $\res(\frac{f_1^{p_1-1}\cdots f_m^{p_m-1}\bullet}{f_1^{p_1}\cdots f_m^{p_m}})=\res(\frac{\bullet}{f_1\cdots f_m})$, which in current language reads 
\begin{equation}\label{ganger}
f_1^{p_1-1}\cdots f_m^{p_m-1}R_{CH}(f^p)=R_{CH}(f).
\end{equation}
It follows that $R^p(f)$ is independent of $p$ if $f$ is a regular sequence. 
In general, however, $R^p(f)$ depends in an essential way on $p$; the set of non-vanishing entries as well as $\ann R^p(f)$ depend on $p$, see Sections ~\ref{monomialcase} and ~\ref{bdiskussion}. Proposition ~\ref{monomB} asserts that if $f$ is monomial, then $R^p$ is independent of $p$ if and only if $f$ is regular. This motivates the following question. 

\begin{questionB}
Suppose that $f=(f_1,\ldots, f_m)$ is a sequence of germs of holomorphic functions at $0\in\C^n$. Let $R^p(f)$ be the Bochner-Martinelli residue current of weight $p$. Is it true that $R^p(f)$ is independent of $p$ if and only if $f$ is regular?
\end{questionB}
Question ~B could be asked also for $\ann R^p(f)$: is it true that $\ann R^p(f)$ is independent if and only if $f$ is a regular sequence?

Lemma 1.2 in ~\cite{BY} asserts that 
\begin{equation}\label{multiplicitet}
\sum_{\I=\{i_1, \ldots, i_n\}\subseteq\{1,\ldots, m\}}
R^p_\I(f)\wedge df_{i_n}\wedge\cdots\wedge df_{i_1}/(2\pi i)^n = e^p(f)[0],
\end{equation}
where $e^p(f)$ is a positive number; in fact each term in \eqref{multiplicitet} is a positive current with support at the origin, see Lemma ~\ref{jwlemma}. Andersson ~\cite{A2} showed that $e^{(1,\ldots, 1)}(f)$ is the Hilbert-Samuel multiplicity of the ideal $\a (f)$. In general $e^p(f)$ depends on $p$, see Example ~\ref{multex}, but it can also happen that $e^p(f)$ is independent of $p$ even if $\ann R^p(f)$ and $R^p(f)$ vary with $p$, as shown in Example ~\ref{konf}.

In general it is hard to compute $R^p(f)$, as well as $\ann R^p(f)$ and $e^p(f)$. However if the $f_j$ are monomials we can give an explicit description of $R^p(f)$ based on ~\cite[Thm.~3.1]{W}. For $A=\{a^1, \ldots, a^m\}\subseteq \Z^n$, let $z^A$ denote the sequence of monomials $z^{a^1},\ldots, z^{a^m}$, where $z^{a^j}=z_1^{a^j_1}\cdots z_n^{a^j_n}$ if $a^j=(a^j_1,\ldots, a^j_n)$. Moreover, for $p\in \mathbb N^m$, let $pA$ denote the set $p A=\{p_1a^1, \ldots, p_ma^m\}\subseteq \Z^n$. 
Given a holomorphic function $g$ we will use the notation $\dbar[1/g]$ for the value at $\lambda=0$ of $\dbar|g|^{2\lambda}$ and analogously by $[1/g]$ we will mean $|g|^{2\lambda}/g|_{\lambda=0}$, that is the principal value of $1/g$.

\begin{thmC}
Suppose that $z^A$ is a sequence of germs of holomorphic monomials at $0\in\C^n$, such that $V(z^A)=\{0\}$.  Let $R^p({z^A})$ be the corresponding Bochner-Martinelli residue current of weight $p$. 
Then 
\begin{equation}\label{formen}
R^p_\I(z^A)= \sgn (A_\I) 
C_\I ~\dbar \left [\frac{1}{z_1^{\alpha^\I_1}}\right ]\wedge\cdots\wedge 
\dbar \left [\frac{1}{z_n^{\alpha^\I_n}}\right ];
\end{equation}
here $\sgn (A_\I)$ is the sign of the determinant of the matrix with rows $a^{i_1},\ldots, a^{i_n}$, $C_\I\geq 0$ is strictly positive if and only if $\I$ is $p$-essential, and $(\alpha^\I_1,\ldots, \alpha^\I_n)=\alpha^\I=\sum_{j\in\I}a^j$.
\end{thmC}
In particular, Theorem ~C implies that 
\begin{equation*}
\ann R^p(z^A)=\bigcap_{\I \text{ p-essential}} (z_1^{\alpha^\I_1},\ldots, z_n^{\alpha^\I_n}).
\end{equation*}

\smallskip

In Section ~\ref{background} we provide some background on Rees valuations, whereas 
the proof of Theorem ~A occupies Section ~\ref{proofs}.  In Section ~\ref{monomialcase} we focus on when $f$ is monomial; we prove Theorem ~C and compute the coefficients $C_\I$ in some special cases. Finally, in Section ~\ref{bdiskussion} we discuss Question ~B and some related questions. 

\smallskip

\noindent
\textbf{Acknowledgment:} I am grateful to Irena Swanson for helpful discussions.

\section{Rees valuations and essential multi-indices}\label{background}
Let $f=(f_1,\ldots, f_m)$ be a sequence of germs of holomorphic functions at $0\in\C^n$, such that $V(f)=\{0\}$.  
The \emph{Rees valuations} of $\a(f)$ are defined in terms of the normalized blowup $\nu:X^+\to(\C^n,0)$ of $\a(f)$, see ~\cite[Ch.II.7]{Ha}. Since $V(\a)=\{0\}$, $\nu$ is an isomorphism outside $0\in\C^n$ and $\nu^{-1}(0)$ is the union of finitely many prime divisors $E\subseteq X^+$. The Rees valuations  of $\a(f)$ are then associated divisorial valuations $\ord_E$ on $\O_0^n$: $\ord_E(g)$ is the order of vanishing of $g$ along $E$.

Let $\pi:X\to (\C^n,0)$ be a log-resolution of $\a(f)$, see ~\cite[Def.~9.1.12]{Laz}. Then, in fact, a divisorial valuation $\ord_E$ is a Rees valuation of $\a(f)$ if and only if the image of the prime divisor $E\subseteq \pi^{-1}(0)$ under the mapping $\Psi=[f_1\circ \pi: \ldots : f_m\circ \pi]: X\dashrightarrow \C\P^{m-1}$ is of (maximal) dimension $n-1$, see ~\cite[p.~332]{Teissier}.

Consider a multi-index $\I=\{i_1,\dots,i_n\}\subseteq\{1,\dots,m\}$. Let $\pi_\I:\C\P^{m-1}\setminus W_\I \to\C\P^{n-1}$, where $W_\I:=\{w_{i_1}=\ldots = w_{i_n}=0\}\subseteq \C\P^n$, be the projection $[w_1:\ldots: w_m]\mapsto [w_{i_1}:\ldots : w_{i_n}]$. 
We say that $\I$ is \emph{essential} with respect to $E$ (and the sequence $f$) if $\Psi(E)\not\subseteq W_\I$ and 
the mapping $\pi_\I\circ \Psi:  E\dashrightarrow \C\P^{n-1}$ is surjective;  
in particular, 
$\ord_E(f_{i_1})=\ldots=\ord_E(f_{i_n})=\ord_E(\a)$. 
Moreover we say that $\I$ is \emph{essential} (with respect to $f$) if $\I$ is essential with respect to at least one exceptional prime. 
Furthermore we say that $\I$ is $p$-\emph{essential} with respect to $E$ (and $f$) if $\I$ is essential with respect to the divisor $E$ and the sequence $f^p$, and that $\I$ is $p$-\emph{essential} (with respect to $f$) if $\I$ is essential with respect to the sequence $f^p$. 

Observe, that if $\I$ is $p$-essential with respect to $E$, then $\vE$ must be a Rees valuation of $\a(f^p)$. Conversely, if $\vE$ is a Rees valuation of $\a(f^p)$, then there exists at least one multi-index ~ $\I$, which is $p$-essential with respect to $E$. However, note that $\I$ can be $p$-essential with respect to more than one divisor $E$, and conversely that there can be several multi-indices $p$-essential with respect to a given $E$. 

Recall that the integral closure of $\a\subseteq\O^n_0$ can be defined in terms of the Rees valuations of $\a$. Indeed, $h\in\O_0^n$ is in $\overline{\a}$ if and only if $\vE(h)\geq \vE(\a)$ for all Rees valuations $\vE$ of $\a$, see for example ~\cite[Ex.~9.6.8]{Laz}.

Given a sequence $f$ and a multi-index $\I=\{i_1,\ldots, i_n\}$, let $f_\I$ denote the sequence $(f_{i_1},\ldots, f_{i_n})$.

\section{Proof of Theorem A}\label{proofs}
The proof of Theorem ~A is very much inspired by and based on (the proofs of) Theorems ~A and ~B in ~\cite{JW} and it also uses Andersson's construction of residue currents in ~\cite{A}.
The following result is Theorem ~B and Lemma 4.3 in ~\cite{JW}.
\begin{lma}\label{jwlemma}
$R_\I(f)\not\equiv 0$ if and only if $\I$ is essential with respect to $f$. Moreover $R_\I(f)\wedge df_{i_n}\wedge\cdots\wedge df_{i_1}/(2\pi i)^n$ is a positive current and its mass is strictly positive if and only if $\I$ is essential.
\end{lma}

\smallskip
We first prove that $R^p_\I(f)\not\equiv 0$ precisely if $\I$ is $p$-essential. 
If  $\I$ is not $p$-essential, then $R_\I(f^p)=0$ by Lemma ~\ref{jwlemma}, and hence in light of \eqref{enkel} $R^p_\I(f)=0$. 
For the converse, note that 
\begin{equation}\label{relation}
R_\I^p(f)\wedge df_{i_n}\wedge\cdots\wedge df_{i_1} = 
\frac{1}{p_{i_1}\cdots p_{i_n}} R_\I(f^p)\wedge df_{i_n}^{p_{i_n}}\wedge\cdots\wedge df_{i_1}^{p_{i_1}}
\end{equation}
by \eqref{enkel}. Lemma ~\ref{jwlemma} asserts that the right hand side of \eqref{relation} is non-vanishing if $\I$ is essential with respect to $f^p$. Thus $R^p_\I(f)\not\equiv 0$ if $\I$ is $p$-essential. 

\smallskip

The inclusion $\ann R^p(f)\subseteq \a (f)$ follows from Andersson's construction of global Bochner-Martinelli residue currents based on the Koszul complex in ~\cite{A}. We  provide (a sketch of) a proof for completeness. 

We identify the sequence $f=(f_1,\ldots, f_m)$ with a holomorphic section of the dual bundle $V^*$ of a trivial vector bundle $V$ over some neighborhood $\mathcal U$ of $0\in\C^n$, endowed with the trivial metric. If $\{e_i\}_{i=1}^m$ is a global holomorphic frame for $V$ and $\{e^*_i\}_{i=1}^m$
is the dual frame, we can write $f=\sum_{i=1}^m f_i e_i^*$. Let $s^p$ be the section $s^p=\sum_{i=1}^m \bar f_i |f_i|^{2(p_i-1)}e_i$, and let  
\begin{equation*}
u^p =
\sum_\ell\frac{s^p\wedge(\dbar s^p)^{\ell-1}}{|f^p|^{2\ell}}.
\end{equation*}
Then $u^p$ is a section of $\Lambda(V\oplus T_{0,1}^*(\mathcal U))$ (where $e_j\wedge d\bar z_i=-d\bar z_i\wedge e_j$), that is clearly well defined and smooth outside $V(f)=\{0\}$, and moreover 
$\dbar|f^p|^{2\lambda}\wedge u^p,$ 
has an analytic continuation as a current to $\Re \lambda > -\epsilon$, see ~\cite{A}. 
Note that the $e_{i_n}\wedge\ldots\wedge e_{i_1}$-coefficient of $R(u^p):=\dbar|f|^{2\lambda}\wedge u^p|_{\lambda=0}$ is just the current $R^p_\I(f)$, and thus in particular, $\ann R(u^p)=\ann R^p(f)$. Let $\nabla=\delta_f-\dbar: \Lambda(V\oplus T_{0,1}^*(\mathcal U))\to \Lambda(V\oplus T_{0,1}^*(\mathcal U))$; here $\delta_f$ denotes interior multiplication by $f$. Then clearly $\nabla u^p=1$ outside $V(f)$. In ~\cite{A} it was proved if $u$ is any section of $\Lambda(V\oplus T_{0,1}^*(\mathcal U))$ that is smooth and satisfies $\nabla u=1$ outside $V(f)$, then the corresponding current $R(u):=\dbar|f|^{2\lambda}\wedge u|_{\lambda=0}$ satisfies that $\ann R(u)\subseteq \a (f)$. We conclude that $\ann R^p(f)\subseteq \a (f)$.

\smallskip

Given a sequence of germs $g_1,\ldots, g_n\in\O^n_0$, let $\jac (g)$ denote the Jacobian determinant $\jac(g)=\left |\frac{\partial g_i}{\partial z_j} \right |_{1\leq i,j\leq n}$. Observe that the coefficient of $df_{i_n}\wedge \cdots \wedge df_{i_1}$ is just $\pm \jac (f_\I)$. Thus in light of \eqref{relation} and Lemma ~\ref{jwlemma}, $\jac (f_\I)\in \ann R^p_\I(f)$ if and only if $R^p_\I(f)\equiv 0$. Given this we can show that $\ann R^p(f)=\a (f)$ implies that $\a (f)$ is a complete intersection ideal by following the proof of Theorem ~A in ~\cite[Section~5]{JW}. 

\smallskip

It remains to prove that the right inclusion in \eqref{eqa} is strict when $n \geq 2$. Given a multi-index $\I=\{i_1,\ldots, i_n\}$, let $P(\I)=\sum_{j=1}^n\frac{1}{p_{i_j}}$. Pick two multi-indices $\I$ and $\J$, such that $P(\I)\geq P(\J)$. We claim that then $R^p_\J(f)\wedge df_{i_n}\wedge\cdots\wedge df_{i_1}$ either vanishes or is a pointmass at the origin.

Let $\pi:X\to (\C^n,0)$ be a log-resolution of $\a (f^p)$. Then $R_\J(f^p)$ is the push-forward of a current $\tilde R$ on $X$, which has support on the exceptional primes with respect to whom $\J$ is essential. 
More precisely, $\tilde R$ can be decomposed as $\tilde R=\sum \tilde R_E$, where the sum is over the exceptional primes $E\subseteq X$, such that $\J$ is essential with respect to $E$, and $\tilde R_E$ has support on $E$, see ~\cite[Section~6]{JW}. 

Let $E_1$ be an exceptional prime, such that $\J$ is essential with respect to $E_1$, and choose local coordinates $\sigma$ on $X$, so that $E_1=\{\sigma_1=0\}$. Moreover, for $2\leq j\leq n$, let $E_j=\{\sigma_j=0\}$, and let $a_j=\text{ord}_{E_j}(f^p)$. Then locally, $\tilde R_{E_1}$ is of the form 
$\dbar [1/\sigma_1^{na_1}]\wedge 
 [1/(\sigma_2^{na_2}\cdots \sigma_n^{na_n})]\wedge \beta$,
where $\beta$ is a smooth form. 
Observe that for $1\leq \ell\leq m$, $\pi^* f_\ell^{p_\ell}$ is divisible by $\sigma_j^{a_j}$ and so $\pi^* f_\ell$ is divisible by $\sigma_j^{\lceil a_j/p_\ell\rceil}$. It follows that 
\[
\pi^*(f_{j_1}^{p_{j_1}-1}\cdots f_{j_n}^{p_{j_n}-1}) \tilde R_{E_1}
=\dbar [1/\sigma_1^{b_1} ]\wedge 
 [1/(\sigma_2^{b_2}\cdots \sigma_n^{b_n}) ] \wedge \beta,
\]
where $b_j\leq a_j P(\J)$. ~A computation following ~\cite[p.~11]{JW} yields that 
\[
\pi^*(df_{i_n}\wedge\cdots\wedge df_{i_1})=
\sigma_1^{c_1-1}(\sigma_2^{c_2}\cdots \sigma_n^{c_n}\gamma + \sigma_1 \delta) 
d\sigma_1\wedge\cdots\wedge d\sigma_n,
\]
where $c_j\geq a_j P(\I)$ and $\gamma$ and $\delta$ are holomorphic functions. Since, by assumption, $P(\I)\geq P(\J)$,  
$\pi^*(f_{j_1}^{p_{j_1}-1}\cdots f_{j_n}^{p_{j_n}-1}) \tilde R_{E_1}
\wedge \pi^*(df_{i_n}\wedge\cdots\wedge df_{i_1})$ is of the form 
$\dbar [1/\sigma_1 ] \wedge d\sigma_1 \wedge \tilde \beta
= 2\pi i [E_1]\wedge \tilde \beta$, where $\tilde\beta$ is a smooth form. Hence 
\begin{multline*}
R^p_\J(f)\wedge df_{i_n}\wedge\cdots\wedge df_{i_1}=\\
\sum_E \pi_* \left (\pi^*(f_{j_1}^{p_{j_1}-1}\cdots f_{j_n}^{p_{j_n}-1}) \tilde R_{E_1} \wedge \pi^*(df_{i_n}\wedge\cdots\wedge df_{i_1}) \right)
\end{multline*}
is a number (possibly $0$) times the Dirac measure at $0$ and the claim is proved.

Now pick a $p$-essential multi-index $\I$, for which $P(\I)=\max_{\J ~~~~~~ p\text{-essential}} P(\J)$. 
Then the non-vanishing entries of $R^p(f)\wedge df_{i_n}\wedge\cdots\wedge df_{i_1}$ are just pointmasses at the origin; in particular,  
$\jac(f_{\I})\m\subseteq \ann R^p(f)$, where $\m$ denotes the maximal ideal in $\O^n_0$. 
Let $E$ be an exceptional prime, such that $\I$ is $p$-essential with respect to $E$. A direct computation gives that 
$\vE(df_{i_1}^{p_{i_1}}\wedge\ldots\wedge df_{i_n}^{p_{i_n}})=n~\vE(f^p)-1$ and $\vE(dz_1\wedge\ldots\wedge dz_n)\geq \sum_{i=1}^n\vE(z_i)-1$. Note that $\vE(z_k)\geq 1$ for $1\leq k\leq n$. 
Since $df_{i_1}^{p_{i_1}}\wedge\cdots\wedge df_{i_n}^{p_{i_n}}=p_{i_1}\cdots p_{i_n} f_{i_1}^{p_{i_1}-1}\cdots f_{i_n}^{p_{i_n}-1}\jac(f_{\mathcal I})dz_1\wedge\cdots\wedge dz_n$ it follows that
\begin{equation*}
\vE(z_kf_{i_1}^{p_{i_1}-1}\cdots f_{i_n}^{p_{i_n}-1}\jac(f_\I))\leq n~\vE(f^p)-n+1=\vE(\overline{(f^p)^n})-n+1
\end{equation*}
for $1\leq k\leq n$. Here we have used that $\overline{\a}$ is the set of all $h\in\O_0^n$, that satisfy $\vE(h)\geq \vE(\a)$ for all Rees valuations $\vE$ of $\a$, see Section ~\ref{background}. 
Hence, if $n\geq 2$, there are elements, for example $z_k \jac(f_\I)$, in $\m \jac (f_\I)$ that are not in $\overline{(f^p)^n}\colon (f_{i_1}^{p_{i_1}-1}\cdots f_{i_n}^{p_{i_n}-1})$. This proves that the first inclusion in \eqref{eqa} is strict and concludes the proof of Theorem ~A.

\section{The monomial case}\label{monomialcase}
Let $z^A=(z^{a^1},\ldots, z^{a^m})$ be a sequence of germs of monomials in $\O^n_0$. Recall that the \emph{Newton polyhedron} $\np(A)$ is defined as the convex hull in $\R^n$ of the set of exponents of monomials in $\a (z^A)$. The Rees valuations of $\a (z^A)$ are monomial and in 1-1 correspondence with the compact facets (faces of maximal dimension) of $\np(A)$. More precisely, the facet $\tau$ with normal vector $\rho=(\rho_1,\ldots, \rho_n)$ corresponds to the monomial valuation $\text{\ord}_\tau(z_1^{a_1}\cdots z_n^{a_n})=\rho_1a_1+\cdots +\rho_na_n$, see for example ~\cite[Thm.~10.3.5]{HS}. 
Given a multi-index $\I$, let $A_\I$ denote the set $\{a^{i_1},\ldots, a^{i_n}\}\subseteq A$ so that $z^{A_\I}$ is the sequence $z^{a^{i_1}}, \ldots, z^{a^{i_n}}$. Moreover, let $\det (A_\I)$ denote the determinant of the matrix with rows $a^{i_1},\ldots, a^{i_n}$. 
It follows that a multi-index $\I$ is essential with respect to $E_\tau$ precisely if $A_\I$ is contained in $\tau$ and $\det (A_\I)\neq 0$; here $E_\tau$ denotes the exceptional prime associated with $\tau$. This means that $\I$ is $p$-essential if and only if $A_\I$ is contained in a facet of $\np(pA)$ and $\det (A_\I)\neq 0$. 

Observe that $z^A$ is regular precisely if $m=n$ and $z^{a^j}$ is of the form $z_j^{b_j}$ (possibly after rearranging the variables). 
Moreover, recall that the integral closure of $\a (z^A)$ is the monomial ideal generated by monomials with exponents in $\np(A)$, see for example ~\cite{Teissier}.

Let us illustrate Theorem ~C with some examples.

\begin{ex}\label{basex}
Let $z^A$ be the sequence of monomials 
$z^A=(z^{a^1},\ldots, z^{a^4})=(z_1^5,z_1^4z_2, z_1^2 z_2^2, z_2^3)$. Then $\np(A)$ has just one compact facet and so $\a(z^A)$ has exactly one Rees valuation, which is the monomial valuation $\vE$ given by $\vE(z_1^{b_1}z_2^{b_2})=3b_1+5b_2$. Moreover the only essential multi-index with respect to $z^{A}$ is $\{1,4\}$ and so Theorem ~C asserts that $R(z^A)=R^{p}(z^A)$, where $p={(1,1,1,1)}$, has one non-vanishing entry and moreover $\ann R(z^A)=(z_1^5,z_2^3)$.

\begin{figure}\label{newtonfigur}
\begin{center}
\includegraphics{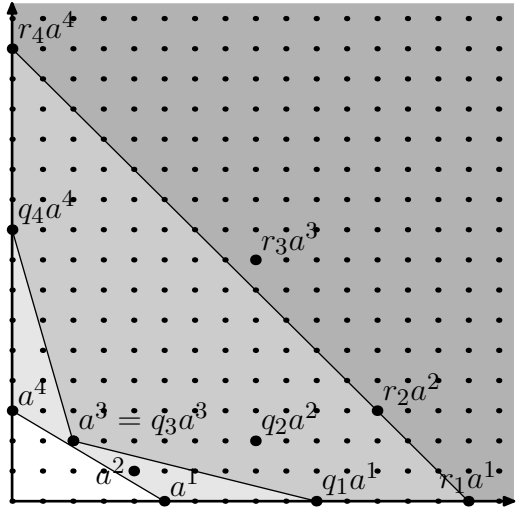}
\caption{The Newton polytopes of the sequences $z^A$ (light grey), $z^{qA}$ (medium grey), and $z^{rA}$ (dark grey) in Example ~\ref{basex}.}
\end{center}
\end{figure}

Let $q=(2,2,1,3)$. Then $\np(qA)$ have two facets, so that $\a(z^{qA})=(z_1^{10}, z_1^8z_2^2, z_1^2z_2^2, z_2^{9})$ has two Rees valuations: $\text{ord}_{E_1}(z_1^{b_1}z_2^{b_2})=b_1+4b_2$ and $\text{ord}_{E_2}(z_1^{b_1}z_2^{b_2})=7b_1+2b_2$. Moreover there are two $q$-essential multi-indices, $\{1,3\}$ and $\{3,4\}$, corresponding to $E_1$ and $E_2$, respectively. It follows from Theorem ~C that
$\ann R^{q}(z^A)=
(z_1^7, z_2^2)\cap (z_1^2, z_2^5)=(z_1^7, z_1^2z_2^2, z_2^5)$.
Note that $\ann R^{p}\not\subseteq \ann R^{q}$ and $\ann R^{q}\not\subseteq \ann R^{p}$, which illustrates that in general no relation between weights $p$ and $q$ are reflected in the relation between $\ann R^p(z^A)$ and $\ann R^q(z^A)$. One can check that by varying the weight $p$ one gets all together 9 different annihilator ideals. Let us consider one more example. 
Let $r=(3,3,4,5)$. Then $\np(rA)$ has one facet, so that $\a(z^{rA})$ has one Rees valuation. However, there are three $r$-essential multi-indices, $\{1,2\}$, $\{1,4\}$, and $\{2,4\}$, and 
$\ann R^{r}(z^A)=(z_1^9,z_2)\cap (z_1^5,z_2^3)\cap (z_1^4,z_2^4)$. In Figure 1 we have drawn $\np(pA)$ and also marked the points in $pA$, for the weights $p$, $q$, and $r$. 

Note that $z^{(q_1-1) a^1}z^{(q_3-1) a^3}=z_1^5$ and $z^{(q_3-1) a^3}z^{(q_4-1) a^4}=z_2^{10}$. It follows that for the weight $q$ the leftmost ideal in \eqref{eqa} is given by
$(z_1^{15}, z_1^{11} z_2, z_1^7z_2^2, z_1^3z_2^3, z_1^2z_2^5, z_1z_2^9, z_2^{12})$ 
and so ones sees directly that the left inclusion in \eqref{eqa} is strict in this case. Note that $\overline{\a (z^A)^2}\not\subseteq \ann R^q(z^A)$, which shows that it is not true in general that $\overline{\a(f)^n}\subseteq \ann R^p(f)$. In Figure 2 the three ideals in \eqref{eqa} are depicted for weights $p$, $q$, and $r$. Note that $\ann R^p(f)$ is strictly included in $\a(f)$ is all three cases.

\begin{figure}\label{annihilatorfigur}
\begin{center}
\includegraphics{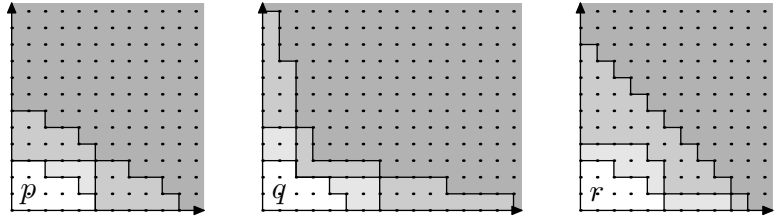}
\caption{The exponent sets of the ideals $\a (z^A)$ (light grey), $\ann R^p(z^A)$ (medium grey) and $\bigcap \overline{(z^A)^2}:(z^{(p_{i_1}-1)a^{i_1}} z^{(p_{i_2}-1)a^{i_2}})$ (dark grey) for weights $p$, $q$, and $r$ in Example ~\ref{basex}.}
\end{center}
\end{figure}

\end{ex}

\begin{ex}\label{tupel}
Let $z^A=(z, z^2)$. Then $\a (z^A)$ is just the maximal ideal $\m \subseteq \O_0^1$. Note that since $n=1$ there is a unique Rees valuation associated with $\a (z^A)$, namely the order of vanishing at the origin. 
For $j\in\mathbb N$, let $p^j=(j,1)$. Then $R(z^A)=R^{p^1}(z^A)=\left (\dbar\left[1/z\right], 0\right)$, $R^{p^2}(z^A)=\left (\dbar\left[1/z\right], \dbar\left[1/z^2\right ] \right)$, and $R^{p^j}(z^A)=\left (0, \dbar\left[1/z^2\right ]\right)$ for $j \geq 3$. It follows that $\ann R=\m$, whereas $\ann R^{p^j}=\m^2$ for $j\geq 2$. 
\end{ex}

Example ~\ref{tupel} shows that in general $R^p(f)$, as well as $\ann R^p(f)$, depends in an essential way on the particular sequence $f$ and not only on the ideal $\a (f)$. Theorem ~A in ~\cite{JW} asserts that $\ann R(f)=\a (f)$ if and only if $\a (f)$ is a complete intersection ideal. Theorem ~A sais that the only if-direction of this statement holds for any $p$, whereas Example ~\ref{tupel} shows that the if-direction fails in general. Moreover, in the monomial case $R(f)$ only depends on $\a(f)$ and not on the particular sequence $f$. Question D in ~\cite{JW} ask whether it is always true (as long as $V(f)=\{0\}$) that $\ann R(f)$ only depends on $\a(f)$.

\subsection{Proof of Theorem C}\label{cbevis}
Theorem 3.1 in ~\cite{W} states that if $\I$ is essential with respect to $z^A$, then $R_\I(f)$ is of the form \eqref{formen}, where $C_\I$ is a nonzero constant. Thus, using \eqref{enkel} and \eqref{ganger}, we conclude that the entries of $R^p(f)$ are of the form \eqref{formen}. 

Assume that $\I$ is $p$-essential. Then by Lemma ~\ref{jwlemma}, \eqref{relation} times $1/(2\pi i)^n$ has strictly positive mass. Note that $dz^{a^{i_n}}\wedge\cdots\wedge dz^{a^{i_1}}=\det(A_\I) dz_n\wedge\cdots\wedge dz_1$. It follows that the left hand side of \eqref{relation} is equal to $(2\pi i)^n C_\I|\det (A_\I)|$, and so $C_\I\geq 0$.

\subsection{The coefficients $C_\I$}\label{coefficients}
Given a sequence of monomials $z^A$ one can find a log-resolution $X_A\to(\C^n,0)$ of $\a (z^A)$, where $X_A$ is a toric variety constructed from the (normal fan of) $\np(A)$, see ~\cite[p.~82]{BGVY}. In ~\cite{W} we computed $R(z^A)$ as the push-forward of a certain current on $X_A$. 
Assume that $\I$ is essential with respect to $E_\tau$, where $\tau$ is a facet of $\np(A)$. According to ~\cite[p.~381]{W}, the coefficient $C_\I$ is of the form $C_\I=\pm \frac{1}{(2\pi i)^{n-1}}  (n-1)! D I$, where $I$ is an integral of the form
\begin{equation*}
I=\int
\frac{\prod_{j=1}^{n-1}|t_j|^{2(c_{j1}+\ldots + c_{jn}-1)}}
{\sum_{k=1}^\ell
\prod_{j=1}^{n-1}|t_j|^{2c_{jk}}} ~d\bar t_1\wedge \cdots \wedge d\bar t_{n-1}\wedge d t_{n-1}\wedge\cdots\wedge d t_1, 
\end{equation*}
for some $n\leq \ell\leq m$ and $\{c_{jk}\}_{1\leq j \leq n-1, 1\leq k\leq \ell}$, and $D$ is the determinant of the matrix with entries $\{d_{jk}\}_{1\leq j,k\leq n}$, where $d_{jk}=c_{jk}$ if $j\leq n-1$ and $d_{nk}=1$. 
The terms in the denominator correspond to the $a^j\in A$ that lie in $\tau$; in particular, $C_\I$ depends only on $\tau\cap A$.  (Assuming that $\I=\{1,\ldots, n\}$ and that $\{a^1,\ldots, a^\ell\}$ are the exponents in $\tau$, then, in the terminology of ~\cite{W}, $c_{jk}=\rho_j\cdot(b_k-a_0)$.) 
In general the integral $I$ is hard to compute; compare to \eqref{integralen2}.

Assume that $\ell=n$ and that $c_{jk}=0$ unless $j=k$, possibly after rearranging the variables $t_j$. Then 
\begin{equation*}
I= \int
\frac{\prod_{j=1}^{n-1}|t_j|^{2(c_{j}-1)}}
{(1+\sum_{j=1}^{n-1}|t_j|^{2c_j})^n}
~d\bar t_1\wedge \cdots \wedge d\bar t_{n-1}\wedge d t_{n-1}\wedge\cdots\wedge d t_1,
\end{equation*}
where $c_j$ just denotes $c_{jj}$. A direct computation gives that 
\begin{equation*}
\int\frac{|s|^{2(N-1)}}{(1+|s|^{2N})^p}d\bar s \wedge ds = 2\pi i \frac{1}{p-1}\frac{1}{N},
\end{equation*}
which implies that $I= \frac{(2\pi i)^{n-1}}{(n-1)!}\frac{1}{c_1\cdots c_{n-1}}$. 
Moreover $D=c_1\cdots c_{n-1}$, and since $C_\I\geq 0$, we conclude that $C_\I=1$. 

The assumption that $\ell=n$ is satisfied precisely if $\I$ is the unique multi-index that is essential with respect to a certain Rees valuation. The assumption that $c_{jk}=0$ for $j\neq k$ is for example satisfied if the normal fan of $\np(A)$ is regular, see ~\cite{Fulton}. It is also satisfied if $n=2$.

Given a facet $\tau$ of $\np(A)$, let $\det(\tau)$ be the normalized volume, that is, $n!$ times the Euclidean volume, of the convex hull of $\tau$ and the origin in $\R^n$. If $\tau$ is simplicial with vertices $b^1,\ldots, b^n$, then $\det(\tau)$ is just (the absolute value of) the determinant of the matrix with rows $b^1,\ldots, b^n$. 
For $n=2$ we have the following description of the coefficients $C_\I$:  
\begin{equation}\label{coffe}
\sum_{A_\I\subseteq\tau} |\det(A_\I)|C_\I=\det(\tau).
\end{equation}
To prove this, recall that if $V(z^A)=\{0\}$, then the Hilbert-Samuel multiplicity $e(z^A)$ of $\a (z^A)$ equals the normalized volume $\vol(\R^n_+\setminus \np(A))$ of the complement in $\R^n_+$ of $\np(A)$, see for example ~\cite{T2}. Observe that $\vol(\R^n_+\setminus \np(A))=\sum \det (\tau)$, where the sum runs over the facets $\tau$ of $\np(A)$. Now \eqref{coffe} follows in light of \eqref{multiplicitet} and the fact that if $\I$ is essential with respect to $E_\tau$, then $C_\I$ depends only on $a^j\in A\cap \tau$. 
\begin{question}
Does \eqref{coffe} hold also when $n>2$?
\end{question}

\begin{ex}\label{multex}
Let $z^A$ and $p$, $q$, and $r$ be as in Example ~\ref{basex}, and let $s=(2,1,1,2)$. From ~\cite{A2} we know that $e^p(z^A)$ is the Hilbert-Samuel multiplicity of $\a (z^A)$. Since there is only one essential multi-index with respect to $z^A$ we can also compute this directly from \eqref{coffe}. Indeed $C_{\{1,4\}}=1$ and so $e^p(z^A)=|\det(A_{\{1,4\}})|=15$. 

Moreover, recall that $\a (z^{qA})$ has two Rees valuations 
and that there is one $q$-essential multi-index associated with each divisor: $\{1,3\}$ and $\{3,4\}$. 
Hence $C_{\{1,3\}}=C_{\{3,4\}}=1$ and so $e^q(z^A)= |\det (A_{\{1,3\}})|+|\det (A_{\{3,4\}})|= 10+6=16$,  that is, the normalized area of the convex hull of $a^1=(5,0)$, $a^3=(2,2)$, and $a^4=(0,3)$. 
Similarly $\a (z^{sA})$ has three Rees valuations and there is one $s$-essential multi-index for each valuation; 
it follows that $e^s(z^A)=17$, see Figure 3. 

Finally $\a (z^{rA})$ has one Rees valuation, but there are three $r$-essential multi-indices. From \eqref{coffe} we know that $C_{\{1,2\}} |\det (A_{\{1,2\}})|+C_{\{1,4\}}|\det (A_{\{1,4\}})|+C_{\{2,4\}}|\det (A_{\{2,4\}})|=|\det (A_{\{1,4\}})|$, which means $C_{\{1,2\}}+ 5C_{\{1,4\}}+ 4C_{\{2,4\}}=5$. However, we cannot say more; in particular, we cannot determine $e^r(z^A)$. 

\begin{figure}\label{hilbertfigur}
\begin{center}
\includegraphics{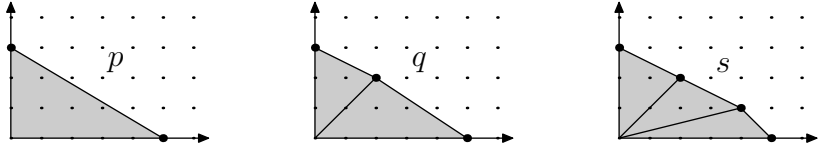}
\caption{The multiplicities $e^p(z^A)$, $e^q(z^A)$, and $e^s(z^A)$ in Example ~\ref{multex}.}
\end{center}
\end{figure}
\end{ex}

\section{Discussion of Question B}\label{bdiskussion}
Theorem ~C allows us to give an affirmative answer to Question ~B in the monomial case. 
Recall that if $\a (f)$ is a complete intersection ideal, then $\a(f)$ is, in fact, generated by $n$ of the $f_j$. This follows for example by Nakayama's Lemma.

\begin{prop}\label{monomB}
Suppose that $z^A=(z^{a^j})_{j=1}^m$ is a sequence of holomorphic monomials at $0\in\C^n$, such that $V(z^A)=\{0\}$.
 Then $R^p(z^A)$ is independent of $p$ if and only if $z^A$ is a regular sequence. 

Moreover, $\ann R^p(z^A)$ is independent of $p$ if and only if for any $\I=\{i_1,\ldots, i_n\}\subseteq\{1,\ldots, m\}$, either $z^{A_\I}$ generates $\a (z^A)$ or $\det (A_\I)=0$. 
\end{prop}
Note that the condition that either $z^{A_\I}$ generates $\a (z^A)$ or $\det (A_\I)=0$ is equivalent to that $\a (z^A)$ is a complete intersection ideal, generated by say $z_1^{b_1},\ldots, z_n^{b_n}$, and that moreover, for $1\leq j\leq m$, $z^{a^j}$ is equal to $z_k^{b_k}$ for some $1\leq k\leq n$. 

\begin{proof}
First, note that the if-directions of the statements in Proposition ~\ref{monomB} are trivially satisfied. Thus we need to prove the only if-directions.

Let $\I'$ be a multi-index defined by that $z^{a_{i_j}}$ is of the form $z_j^{b_j}$, where $b_j$ is the smallest number such that $z_j^{b_j}$ is among the entries of $z^A$. Without loss of generality we may assume that $\I'=\{1,\ldots, n\}$.  Choose $p\in \mathbb N^m$, so that $p_{i}=1$ if $i\leq n$ and $p_i >>1$ otherwise. Then $\I'$ is the unique $p$-essential multi-index. 

Assume that $m>n$ and choose $j$, such that $n<j\leq m$. Moreover, choose $q\in \mathbb N^m$ such that $a^j$ lies in the one of the compact facets of the boundary of $\np(qA)$. For example, let $q$ be defined by $q_i=|a^1|+\ldots +|a^{i-1}|+|a^{i+1}|+\ldots +|a^m|$, where $|a^j|=a^j_1+\ldots + a^j_n$. Then $j$ is contained in a $q$-essential multi-index, say $\J$. It follows that $R^q_\J(f)\neq 0$, whereas $R^p_\J(f)=0$. Hence $R^p(z^A)\neq R^q(z^A)$ and we have proved the first part of Proposition ~\ref{monomB}. 

Next, assume that there is an $a^j\in A$ that is not equal to any of $b_1,\ldots, b_n$. Since $V(z^A)=\{0\}$, at least one of the entries of $a^j$ is positive, say $a^j_k>0$. Let $\J=\{1,\ldots, k-1, k+1, \ldots, n,j\}$. Then $\det (A_\J)\neq 0$, which means that we can find a weight $q$ such that $\J$ is $q$-essential; for instance we can take $q$ as above. By assumption, $a^j_k>b_k$ or $a^j_i>0$ for some $i\neq k$. In both cases $\sum_{j\in\J}a^j\neq \sum_{j\in\I}a^j$, and thus $\ann R^q_\J(z^A)\neq\ann R^p_\I(z^A)$, where $p$ is as above. This proves the second part of Proposition ~\ref{monomB}. 
\end{proof}

\smallskip

Observe that a sufficient condition for Question ~B to be true would be that the set of $p$-essential multi-indices is independent of $p$ if and only if $f$ is a regular sequence. As we saw in the above proof this is true if $f$ is monomial, but we do not know if it holds in general. When $f$ is monomial, the essential multi-indices are rather special. For example, a multi-index can be essential with respect to at most one Rees valuation, which is not the case in general. Indeed, if $m=n$, then $\I=\{1,\ldots, n\}$ is essential with respect to all Rees valuations (and there can be more than one Rees valuation). The following example illustrates another phenomenon, which does not occur in the monomial case. 

\begin{ex}\label{irena}
Let $f=(z_1^4-z_2^4, z_1^2z_2, z_1z_2^2)$. Then $\a (f)$ has three Rees valuations, namely the monomial valuations $\text{ord}_{E_1}(z_1^{b_1}z_2^{b_2})=b_1+b_2$, $\text{ord}_{E_2}(z_1^{b_1}z_2^{b_2})=2b_1+b_2$, $\text{ord}_{E_3}(z_1^{b_1}z_2^{b_2})=b_1+2b_2$. Note that $\{2,3\}$, $\{1,3\}$ and $\{1,2\}$ are the unique essential multi-indices with respect to $\text{ord}_{E_1}$, $\text{ord}_{E_2}$, and $\text{ord}_{E_3}$, respectively. Note that this situation cannot happen if $f_j$ are all monomials. 

Let $q=(1,2,2)$. Then $\a (f^q)=(z_1^4-z_2^4, z_1^4z_2^2, z_1^2z_2^4)$ has four Rees valuations, $\text{ord}_{E_1}, \ldots, \text{ord}_{E_4}$. To see this, note that after blowing up the origin once, the strict transform of $\a(f^q)$ has support at four points $x_1,\ldots, x_4$. The divisor $E_j$ is obtained by further blowing up $x_j$ twice. A computation yields that $\{1,2\}$ and $\{1,3\}$ are both $q$-essential with respect to $E_j$ for $1\leq j\leq 4$, whereas $\{2,3\}$ is not $q$-essential. Hence $R(f)\neq R^q(f)$. 
\end{ex}

\smallskip

Note that $\det (A_\I)=0$ is equivalent to that $dz^{a^{i_1}}\wedge\cdots\wedge dz^{a^{i_n}}$ vanishes identically, which in turn implies that $\I$ is not $p$-essential for any $p\in\mathbb N^m$, see for example Lemma ~\ref{jwlemma}. 
This motivates the following version of Question ~B.

\begin{questionB'}
Is it true that $\ann R^p(f)$ is independent of $p$ if and only if for any $\I=\{i_1,\ldots, i_n\}$, either $f_\I$ generates $\a (f)$ or the form $df_{i_1}\wedge\cdots\wedge df_{i_n}$ vanishes identically? 
\end{questionB'}
Let us mention some partial answers to Question B'. 
Theorem ~C in ~\cite{JW} asserts that if $\a(f)$ is a complete intersection ideal, then $R_\I(f)$ is a constant times the Coleff-Herrera product $R_{CH}(f_\I)$ if $f_\I$ generates $\a(f)$ and $0$ otherwise. Using this and \eqref{ganger} one can check that $\ann R^p(f)$ is independent of $p$ if $\a(f)$ is a complete intersection ideal, generated by say $f_1,\ldots, f_n$, and moreover for $j>n$, $f_j$ is equal to (a constant times) one of the $f_k$ for $1\leq k\leq n$; compare this to (the discussion right after) Proposition ~\ref{monomB}.  

\begin{ex}\label{olika}
Let $f=(z,w,z+w)$. Then $\a (f)$ is just the maximal ideal in $\O^2_0$, which is clearly a complete intersection ideal, and thus by Theorem ~C in ~\cite{JW}, $\ann R(f)=\a(f)$. Note that any choice of $f_i$ and $f_j$ generate $\a(f)$, so $f$ satisfies the condition in Question ~B'.

Let $p=(3,3,3)$. Observe that $\a (f^p)=(z_1^3,z_2^3, z_1^2z_2+z_1z_2^2)$ is not a complete intersection ideal. A computation yields that
\[
R_{\{1,3\}}(f^p)=A_1\dbar[1/z^5]\wedge \dbar [1/w]+
A_2\dbar [1/z^4]\wedge \dbar [1/w^2]
+A_3\dbar [1/z^3]\wedge \dbar [1/w^3],
\]
for some constants $A_1$, $A_2$, and $A_3$. It follows that 
$R^p_{\{1,3\}}(f)=(A_1+2A_2+A_3)\dbar[1/z]\wedge\dbar[1/w]$.  In fact, also the other entries of $R^p$ are of this form and so $\ann R^p=\a(f)$. 
\end{ex}

Note that if there is a subsequence $f_\I=(f_{i_1},\ldots, f_{i_n})$ of $f$ such that $V(f_\I)=\{0\}$, then by choosing $p_j=1$ if $j\in\I$ and $p_j >>1$ for $j\notin\I$, the only non-vanishing entry of $R^p(f)$ is $R^p_\I(f)$, which is a constant times $R_{CH}(f_\I)$. Thus, given that there exists such an $f_\I$, $R^p(f)$ is not independent of $p$ as soon as, for example, there is another multi-index $\J$, such that $V(f_\J)=\{0\}$, or as soon as $\ann R(f)$ is not a complete intersection ideal. One can, however, not always find such an $f_\I$, as the following example shows. 

\begin{ex}\label{slutexempel}
Let $f=(z_1z_2,z_1(z_1+z_2),z_2(z_1+z_2))$. Then $V(f_{\I})$ is a line through the origin for all $\I=\{i_1,i_2\}$; in particular, $V(f_{\I})\neq \{0\}$. 
Moreover, $\a(f)$ is the (monomial) ideal $\m^2$, where $\m$ is the maximal ideal in $\O_0^2$. Thus the only Rees valuation of $\a(f)$ is the order of vanishing at the origin and so $R(f)$ can be computed by blowing up the origin once. Note that all multi-indices $\I=\{i_1,i_2\}$ are essential. 
Let $R^{\ell, k}$ denote the current $\dbar[1/z_1^\ell]\wedge\dbar[1/z_2^k]$, and let 
\begin{equation}\label{integralen2}
C_j=\frac{1}{2\pi i}\int \frac{|t|^{2j}d\bar t \wedge dt}
{(|t|^2+|1+t|^2+|t(1+t)|)^2}.
\end{equation} 
Then, a computation yields that $R_{\{1,2\}}(f)=-C_0R^{3,1}$, $R_{\{1,3\}}(f)=2C_2R^{1,3}$, and 
$R_{\{2,3\}}(f)=C_0R^{3,1}+2C_1R^{2,2}+C_2R^{1,3}$. It follows that $\ann R(f)=\m^3$. 

Let $p=(2,1,1)$. Then $\a(f^p)$ has two Rees valuations, $\text{ord}_{E_1}$ and $\text{ord}_{E_2}$, where $E_1$ is the exceptional divisor obtained by blowing up the origin once, whereas $E_2$ is obtained by further blowing up a point on $E_1$ twice. Moreover, $\{2,3\}$ is essential with respect to $E_1$ and $\{1,2\}$ and $\{1,3\}$ are essential with respect to $E_2$. A computation gives that $R^p_{\{1,2\}}(f)=R^p_{\{1,3\}}(f)=-1/2(R^{3,1}-R^{2,2}+R^{1,3})$ and $R^p_{\{2,3\}}(f)=A^{3,1}R^{3,1}+A^{2,2}R^{2,2}+A^{1,3}R^{1,3}$, where $A^{i,j}>0$. 

Note that $R^p_\I(f)\neq R_\I(f)$, as well as $\ann R^p_\I(f)\neq \ann R_\I(f)$, for all $\I$. Moreover, note that $\ann R(f)$ is strictly included in $\ann R^p(f)$. Indeed, 
$(A^{2,2}+A^{1,3})z_1^2 + (A^{1,3}-A^{3,1})z_1z_2 - (A^{3,1}+A^{2,2})z_2^2\in\ann R^p(f)\setminus \ann R(f)$. 
\end{ex}

\subsection{Related questions}
Question ~B could be posed also for the currents \eqref{multiplicitet}. The following example shows that $e^p(f)$ does not necessarily vary with $p$ even if $R^p(f)$ and $\ann R^p(f)$ do. 

\begin{ex}\label{konf}
Let $z^A=(z_1^2,z_1z_2, z_2^2)$. Then by varying $p$ there are three different possibilities of $p$-essential multi-indices. First, all three multi-indices $\I=\{i_1,i_2\}$ can be $p$-essential, which for example is the case for $p=(1,1,1)$. Next, for $p=(1,2,1)$, $\{1,3\}$ is the only $p$-essential multi-index, and for $p=(2,1,1)$, the $p$-essential multi-indices are $\{1,2\}$ and $\{2,3\}$. In the first situation, by \cite{A2}, $e^p(z^A)$ is the Hilbert-Samuel multiplicity of $\a(z^A)$, which is equal to $\vol (\R^n_+\setminus \np(A))=|\det(A_{\{1,3\}})|=4$. In light of \eqref{coffe} it is not hard to check that this is holds true also if $p$ is another weight such that all $\I$ are $p$-essential. In the latter two cases, by Section ~\ref{coefficients}, the coefficients $C_\I$ are all $1$, when $\I$ is $p$-essential.  It follows that $e^p=|\det (A_{\{1,3\}})|=4$ and $e^p=|\det (A_{\{1,3\}})|+|\det (A_{\{2,3\}})|=2+2$, respectively, so in fact $e^p(z_A)$ is independent of $p$. 
\end{ex}

One can also ask in what sense $R^p(f)$ and $\ann R^p(f)$ depend on $p$, once $\I$ is $p$-essential. In the monomial case $\ann R^p_\I(f)$ is fix as long as $\I$ is essential but the coefficient $C_\I$ in \eqref{formen} vary in general. Indeed, in Example ~\ref{konf} above, for $p=(1,1,1)$, $C_{\I}$ are all strictly between $0$ and $1$, whereas in the latter cases they are either $0$ or $1$. In general, also $\ann R^p_\I(f)$ varies with $p$, see Example ~\ref{slutexempel} above. Computations, such as in Example \ref{slutexempel}, suggest that in general there may be infinitely many different annihilator ideals $\ann R^p_\I(f)$ and $\ann R^p(f)$ as $p$ varies over $\mathbb N^m$. This contrasts the monomial case, where there are always finitely many different ideals $\ann R^p(f)$.

\def\listing#1#2#3{{\sc #1}:\ {\it #2},\ #3.}

\end{document}